\newcommand*{\memfontfamily}{pazo}
\newcommand*{\memfontpack}{mathpazo}

\documentclass[10pt,
onecolumn,oneside,a5paper,article
,frenchb,italian,german,swedish,latin%
,british%
]{memoir}

\usepackage[mathpazo
,eugreek
,optima
,itsans
,lmodern
]{memmana2}

  \hypersetup{
breaklinks=true,%
    pdftitle = {Conjectures and questions in convex geometry},
    pdfauthor = {PierGianLuca  Porta Mana}
  }

\usepackage{mathdots}

\renewcommand{\multicitedelim}{\addsemicolon\space}
\bibliography{manabibliography}
\renewcommand{\bibfont}{\small}

\DeclareBibliographyCategory{1}
\newcommand{\citep}{\parencites}
\newcommand{\citey}{\parencites*}
\renewcommand{\cite}{\citep}

\setlxvchars
\setxlvchars
\settypeblocksize{*}{\lxvchars}{1.618}
\setulmargins{*}{*}{*}
\setlrmargins{*}{*}{*}
\setheaderspaces{*}{*}{*}
\checkandfixthelayout

\providecommand{\firstname}[1]{#1}
\providecommand{\surname}[1]{#1}
\providecommand{\affiliation}[1]{{\footnotesize#1}}
\providecommand{\epost}[1]{{\footnotesize #1}}
\providecommand{\email}[1]{{\textsf{\href{mailto:#1}{#1}}}}
\providecommand{\pacs}[1]{{\footnotesize\textsc{PACS} numbers: #1}}
\providecommand{\msc}[1]{{\footnotesize\textsc{MSC} numbers: #1}}

\ifdraftdoc
\usepackage
{datetime}
\newdateformat{mydate}{\THEDAY\ \monthname[\THEMONTH] \THEYEAR}
\fi

\newenvironment{acknowledgements}{\chapter*{Acknowledgements}\addcontentsline{toc}{section}{Acknowledgements}}{\par}

\chapterstyle{reparticle}
\copypagestyle{manaart}{plain}
\makeheadrule{manaart}{\headwidth}{0.75\normalrulethickness}
\makeoddhead{manaart}{%
\footnotesize\textit{Porta Mana}}{}{
\footnotesize\textit{Conjectures and questions in convex geometry}}
\makeoddfoot{manaart}{}{\thepage}{}
\makeatletter
\newcommand\addprintnote{%
\begin{picture}(0,0)%
\put(0,-14){
\makebox(0,0){
{\tiny%
This document is optimized for on-screen reading and 2-pages-on-1-sheet
printing on A4 or Letter paper}}
}%
\end{picture}%
}
\makeoddfoot{plain}{}{\makebox[0pt]{\thepage}\addprintnote}{
}
\makeoddhead{plain}{}{}{\footnotesize
Report no.\ pi-qf-223
}

\title{Conjectures and questions in convex geometry\\[1.5\jot]\large 
of interest for quantum theory\\ and other physical statistical theories
}
\author{\firstname{P.G.L.}\ \surname{Porta\,Mana}
}
\postauthor{\\[-.5ex]%
\affiliation{Perimeter Institute for Theoretical Physics, Canada} 
\\[-1ex]
\epost{\textless\email{lmana AT pitp.ca}\textgreater
}\end{tabular}\par\end{center}}

\predate{\begin{center}\footnotesize}
\date{16 May 2011
{
\\ (first drafted 21 November 2009)}
}
\postdate{\end{center}}

\usepackage{breakurl} 
\usepackage{pifont}

\theoremstyle{plain}
\newtheorem{theorem}{Theorem}
\theoremstyle{remark}

\newtheorem{quest}{Question}
\theoremstyle{definition}
\newtheorem{example}{Example}
\newtheorem{cexample}[example]{Counter-example}
\newtheorem{definition}{Definition}
\newtheorem{fact}{Fact}

\newtheorem{conj}{Conjecture}
\newtheorem{tconj}{Physical conjecture}

\newcommand{\yC}{\mathcal{C}} 
\newcommand{\yO}[1]{\mathcal{P}_{#1}}
\newcommand{\yOC}{\yO{\yC}}
\newcommand{\yOT}{\yO{\yT}}
\newcommand{\yom}{v}
\newcommand{\yoz}{\yom_0}
\newcommand{\you}{\yom_\text{u}}
\newcommand{\yT}{\varDelta}
\newcommand{\yTT}{\bar{\yT}}
\newcommand{\yOTT}{\yO{\yTT}}

\newcommand{\yf}{F}
\newcommand{\yg}{G}
\newcommand{\yfp}{\pi}
\newcommand{\ygp}{\gamma}
\newcommand{\yff}{\bar{\yf}}
\newcommand{\ygg}{\bar{\yg}}
\newcommand{\ya}{\lambda}
\newcommand{\yga}{\alpha}
\newcommand{\ypa}{\psi_a}
\newcommand{\ypb}{\psi_b}
\newcommand{\yozz}{d_0}
\newcommand{\youu}{d_\text{u}}
\newcommand{\ydr}{\varTheta}
\newcommand{\yer}{\varEpsilon}
\newcommand{\yb}{x}
\newcommand{\yx}{\bm{q}}
\newcommand{\yll}{L}

\newcommand{\yc}{y_x}
\newcommand{\QEM}
{\ding{167}}
\newcommand{\qem}{\leavevmode\unskip\penalty9999 \hbox{}\nobreak\hfill
\quad\hbox{\QEM}}
\DeclareMathOperator{\aff}{aff}
\DeclareMathOperator{\conv}{conv}

\setfloatadjustment{figure}{\footnotesize\centering}

\begin{document}
\selectlanguage{british}
\hyphenation{ 
Pre-sent Pre-sent-ed Pre-sent-ing Pre-sents Li-ce-o Scien-ti-fi-co Sta-ta-le ca-glia-ri Con-tin-u-um
Quan-tum Be-tween Phe-nom-e-non Mac-ro-scop-ic Mi-cro-scop-ic Sub-space
Sub-spaces Mo-men-tum Mo-men-ta Ther-mo-me-chan-ics Ther-mo-me-chan-i-cal
Meso-scop-ic Elec-tro-mag-net-ic Ve-loc-i-ty Ve-loc-i-ties Gal-i-le-an Gal-i-le-ian
}

\firmlists*
\maketitle
\abslabeldelim{:\quad}
\setlength{\abstitleskip}{-\absparindent}
\abstractrunin
\begin{abstract}
  Some conjectures and open problems in convex geometry are presented, and
  their physical origin, meaning, and importance, for quantum theory and
  generic statistical theories, are briefly discussed.
  \\[2\jot]
  \msc{52B11,14R99,81P13}\\
  \pacs{02.40.Ft,03.65.Ta,05.90.+m}
\end{abstract}


\newrefsegment
\selectlanguage{british}

\chapter{Introduction}
\label{cha:intro}

In this note I present a couple of conjectures and open problems in convex
geometry, wishing that they will raise the interest of geometers and be
solved soon.

What is the origin of these conjectures and open problems? There is a
branch of Bayesian probability theory, called the theory of statistical
models, that studies the probabilistic relations among particular sets of
propositions or variables \citep[][and refs
therein]{holevo1976_t1978,holevo1980_t1982,holevo1985,portamana2011c}. Its
range of applications is therefore as vast and diverse as that of
probability theory. One of these applications, which is gaining the
interest of more and more researchers in physics, mathematics, and
statistics, is the study of the probabilistic and communication-theoretic
features of quantum theory and other physical statistical theories, and of
how these features can be mimicked by or emerge from a classical theory.

Convex geometry is one of the main mathematical structures at the core of
the theory of statistical models. So with this theory we can translate some
theorems and conjectures about quantum and non-quantum theories and their
relations with classical theories into theorems and conjectures about convex
geometry, and vice versa. It is these physical conjectures that are here
presented, translated in strictly convex-geometric terms.

The presentation follows standard notation and terminology
\citep{coxeter1948,coxeter1961_r1969,valentine1964,gruenbaum1967_r2003,klee1971,broendsted1983,webster1994,portamana2011}.
Some definitions are presented in the next section, in particular the
notion of a statistical model; then the notion of \emph{refinement} of a
statistical model is presented in \sect~\ref{sec:maps} with many
illustrative examples: this is the notion to which are connected the open
problems and conjectures presented in \sect~\ref{sec:conjectures}. Physical
motivation and meaning are discussed in an appendix.

\chapter{Definitions}
\label{sec:affintro}

In this section we consider a compact convex space $\yC$ of dimension
$n$.

\begin{definition}
  The \emph{convex-form space} $\yOC$ of the convex space $\yC$ is the set
  of all affine forms on $\yC$ with range in $\clcl{0,1}$, called
  \emph{convex forms}:
\begin{equation}
  \label{eq:outcomesp}
  \yOC \defd
\set{\yom \colon \yC \to \clcl{0,1} \st \text{$\yom$ is affine}}.
\end{equation}
Convex combination is easily defined on this set, which is thus a convex
space itself. A vector sum and difference can also be naturally defined
but the set is not closed under them. The forms $\yoz \colon a \mapsto 0$
and $\you \colon a \mapsto 1$ are called \emph{null-form} and
\emph{unit-form}. The action of a form $\yom$ on a point $a$ is denoted by
$\yom \inn a$. We call $\yom$ an \emph{extreme form} if it is an extreme
point of $\yOC$.
\end{definition}
Later on we shall on occasion write vector differences of convex forms when
their result is still a convex form.

\begin{figure}[!t]
  \centering
  \includegraphics[width=.8\columnwidth]{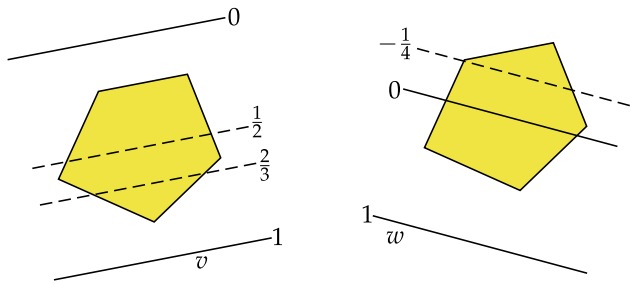}
  \caption{\small The contour hypersurfaces of an affine form are parallel
    hyperplanes; such a form is completely determined by assigning the two
    hyperplanes corresponding to the values $0$ and $1$. To be a convex
    form, these two hyperplanes must not intersect the interior of the
    convex space on which the form is defined. On the left, $v$ is a convex
    form for the pentagonal convex space; two lines are indicated where $v$
    has values $1/2$ and $2/3$; no points of the space yield the values $0$
    or $1$. On the right, $w$ cannot be a convex form (although it is an
    affine form) because it assigns strictly negative values to some points
    of the convex space; this happens because the $0$-value line cuts the
    convex space.}
\label{fig:plausform}
\end{figure}

It is useful to recall that a non-constant convex form $v$ on $\yC$ is
determined by an ordered pair of $(n-1)$-dimensional, parallel hyperplanes
non-intersecting the interior of $\yC$, as explained in
fig.~\ref{fig:plausform}.

%

\begin{fact} If the convex space $\yC$ has dimension $n$ then its
  convex-form space $\yOC$ has dimension $n+1$. If $\yC$ is a polytope, so
  is $\yOC$. The convex structure of $\yOC$ is determined by that of $\yC$,
  and its affine span $\aff\yOC$ is the space of affine forms on $\aff\yC$.
  $\yOC$ is a bi-cone whose vertices are the null-form $\yoz$ and the
  unit-form $\you$, and this bi-cone is centro-symmetric with centre of
  symmetry $(\yoz +\you)/2$. The number of extreme points of $\yOC$ besides
  $\yoz$ and $\you$ is determined by the structure of the faces of $\yC$;
  \eg, if $\yC$ is a two-dimensional polytope, the number of extreme points
  of $\yOC$ equals $2m+2$, where $m$ is the number of bounding directions
  of $\yC$.
\end{fact}

\begin{figure}[!pt]
  \centering
  \includegraphics[width=9.5\columnwidth/10]{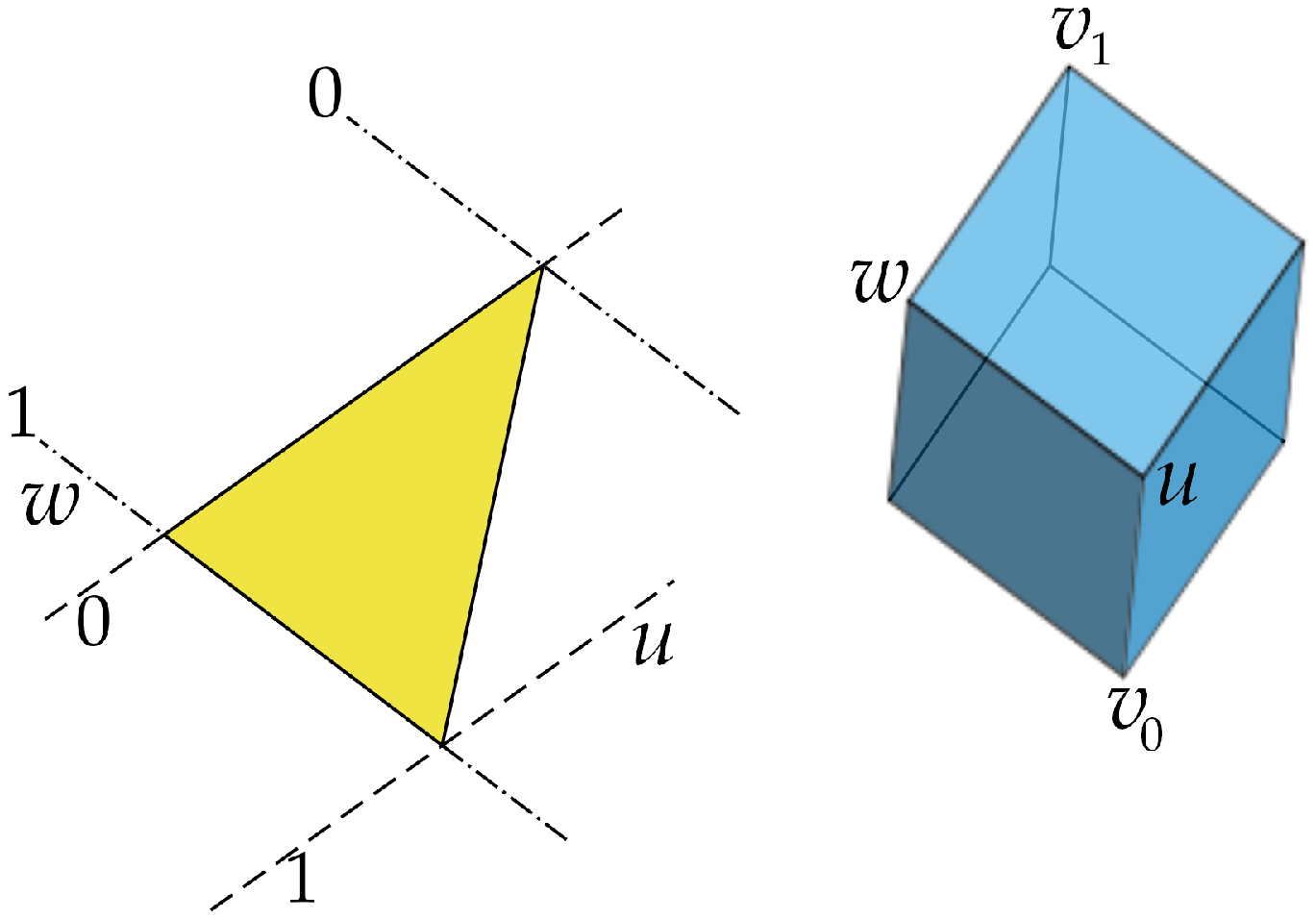}
  \includegraphics[width=9.5\columnwidth/10]{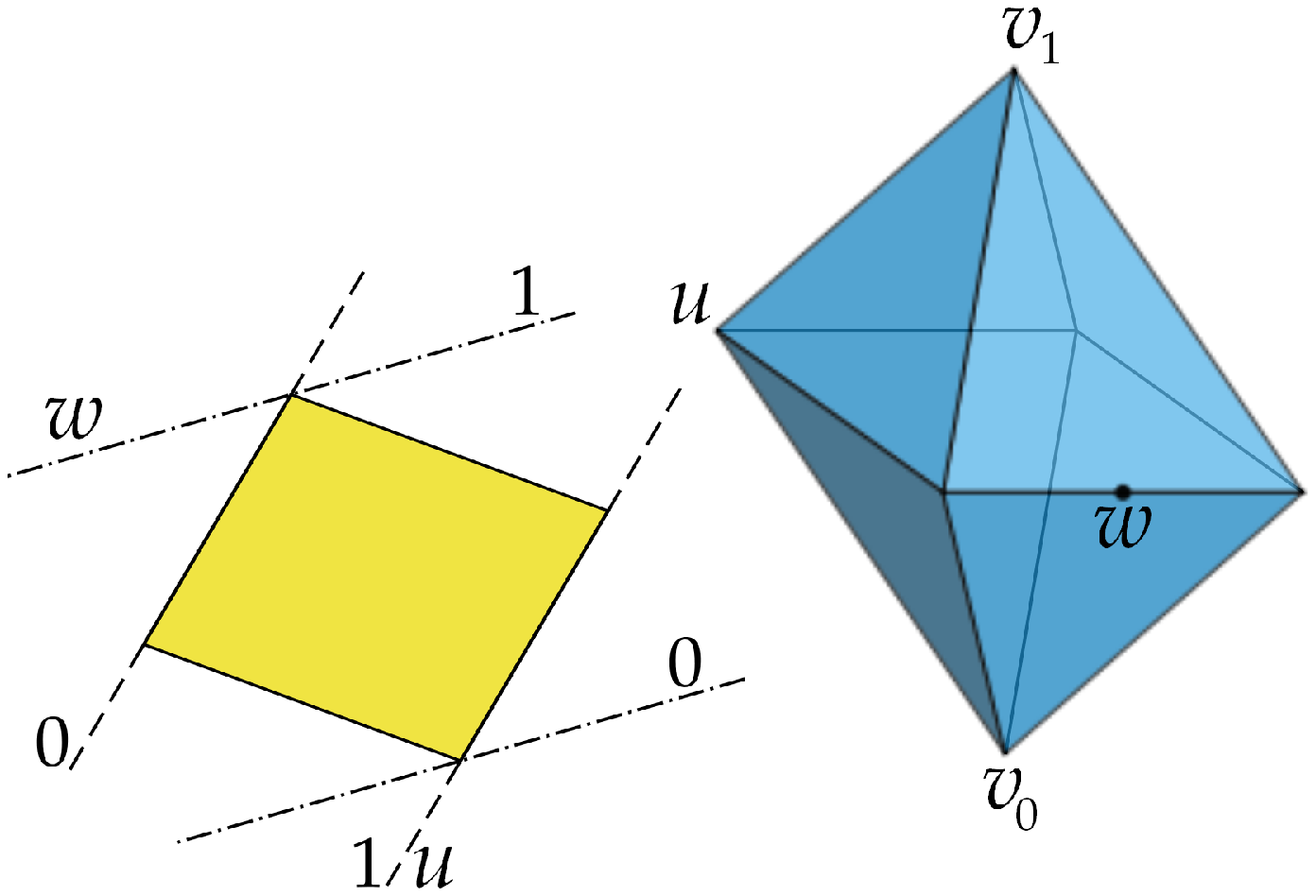}
  \caption{\small Two two-dimensional convex spaces, on the left, with their
    three-dimensional convex-form spaces, on the right. They constitute two
    statistical models. The convex forms $w$, $u$ are represented as pairs
    of parallel lines on the convex spaces and as points on the convex-form
    spaces. $v_0$ and $v_1$ are the null- and unit-forms.}
\label{fig:out_sp}
\end{figure}

For example, the convex-form space of an $n$-dimensional simplex is an
$(n+1)$-dimensional parallelotope (with $2^{n+1}$ extreme points), and that
of a parallelogram is an octahedron, as shown in fig.~\ref{fig:out_sp};
that of a pentagon is a pentagonal trapezohedron, shown in
fig.~\ref{fig:pentag}.

The following definition introduces the most important mathematical objects
of our study:
\begin{definition}
  A \emph{statistical model} is a pair $(\yC, \yOC)$ of a convex polytope
  and its convex-form space; its \emph{dimension} is simply the dimension
  of $\yC$. A \emph{simplicial} statistical model is one in which $\yC$ is
  an $n$-dimensional simplex $\yT$, and therefore $\yOT$ is an
  $(n+1)$-dimensional parallelotope.
\end{definition}
The name, especially the adjective `statistical', is admittedly sibylline,
but it rightly prophesies a connexion with probability theory.

\pagebreak[1]
\chapter{Refinement of a statistical model}
\label{sec:maps}

In this section we consider compact convex \emph{polytopes} of dimension
$n$; \ie, we are assuming that our convex spaces have a finite number of
extreme points.

In convex geometry it is a well-known fact that any polytope can be
obtained as a section or a projection of a usually higher-dimensional
simplex (Gr\"unbaum \citey{gruenbaum1967_r2003}, \sect~5.1, theorems~1 and
2). Both these kinds of `correspondence' $\yC \leftrightsquigarrow \yT$
between a polytope $\yC$ and the simplex $\yT$ from which the polytope is
obtained have the following characteristics:
\begin{enumerate}[C1.]
\item each point of $\yC$ has at least one corresponding point in $\yT$,
  \ie, the correspondence $\yC \leftrightsquigarrow \yT$ is defined on all
  $\yC$;
\item there may be points of $\yT$ with no corresponding points in $\yC$,
  \ie, the correspondence $\yC \leftrightsquigarrow \yT$ needs not be
  defined on all $\yT$;
\item several points in $\yT$ can correspond to one and the same point in
  $\yC$, \ie, $\yC \leftrightsquigarrow \yT$ can be one-to-many;
\item at most one point in $\yC$ can correspond to one in $\yT$, \ie, $\yC
  \leftrightsquigarrow \yT$ cannot be many-to-one;
\item a convex combination of points in $\yT$ corresponds to the same
  convex combination of the corresponding points in $\yC$, when the latter
  are defined, \ie, $\yC \leftrightsquigarrow \yT$ is affine.
\end{enumerate}
These characteristics mathematically pin down the $\yC \leftrightsquigarrow
\yT$ correspondence as a partial, onto, affine map from $\yT$ to $\yC$ (we
denote partial maps by hooked arrows):
\begin{align}
\yf: \yT \xhookrightarrow{\text{onto, affine}} \yC,
\label{eq:generalS}
\end{align}
\begin{itemize}
  \item `onto' = it covers $\yC$, from C1;
  \item `partial' = it needs not be defined on all $\yT$, from C2;
  \item `affine' = if $\yf$ is defined on $a, b$ then $\yf[\ya a + (1-\ya)
    b] = \ya \yf(a) + (1-\ya) \yf(b)$, from C5.
  \item `map' = it is many-to-one or one-to-one, but not one-to-many, from
    C3 and C4.
\end{itemize}

Intuitively, this map says that the simplex $\yT$ has a `finer' structure
than the polytope $\yC$, or that $\yC$ is a `coarser' image of $\yT$
because parts of the latter are either missing or not distinguishable in
$\yC$. We might call $\yT$ a `simplicial refinement' of $\yC$. The
projection or section of a simplex are particular cases of this map.

\bigskip

It is natural to try to generalize this kind of construction and its
associated theorems from polytopes to statistical models: given a
statistical model $(\yC, \yOC)$, one asks whether it can be obtained from a
simplicial statistical model $(\yT, \yOT)$, considered as a `refinement'.

More precisely, we want a correspondence $\yC \leftrightsquigarrow \yT$
between the convex spaces, one $\yOC \leftrightsquigarrow \yOT$ between
their convex-form spaces, and we want both correspondences to satisfy
requirements analogous to C1--C5. Moreover, we clearly want  these
correspondences to preserve the action of convex forms on the respective
convex spaces in the two statistical models.

This means that the two correspondences have to be expressed by partial,
surjective, affine maps from $\yT$ to $\yC$ and from $\yOT$ to $\yOC$
\begin{align}
&\yf: \yT \xhookrightarrow{\text{onto, affine}} \yC,
\label{eq:generalF}
\\
&\yg: \yOT \xhookrightarrow{\text{onto, affine}} \yOC
\label{eq:generalG}
\end{align}
that satisfy
 \begin{multline}   \label{eq:compatib}
   \yg(\yom)\inn \yf(a) = \yom \inn a
\\ \text{for all $a \in \yT$, $\yom \in \yOT$ 
 on which $\yf$, $\yg$ are defined}.
 \end{multline}

We are thus led to the following
\begin{definition}
  A \emph{simplicial refinement} of a statistical model
  $(\yC, \yOC)$ is a set $(\yT, \yOT, \yf, \yg)$ where:
\begin{enumerate}[I.]
\item $(\yT, \yOT)$ is a simplicial statistical model,
\item $\yf\colon \yT \hookrightarrow \yC$ is partial, onto, and affine,
\item $\yg\colon \yOT \hookrightarrow \yOC$ is partial, onto, and affine,
\item  $\yf$ and $\yg$ are such that $\yg(\yom)\inn \yf(a) = \yom \inn
  a$ on their domains of definition.
\end{enumerate}
\end{definition}
We shall often omit the adjective `simplicial' when speaking about a
simplicial refinement.

\begin{example}\label{ex:simplic}

Consider the statistical model where $\yC$ is a parallelogram and $\yOC$ an
octahedron, as at the bottom of fig.~\ref{fig:out_sp}
or~\ref{fig:stat_red}. A simplicial refinement is given by $(\yTT, \yOTT,
\yfp,\ygp)$, where: $\yTT$ is a tetrahedron; $\yOTT$ a hypercube; the map
$\yfp$ is the projection of the tetrahedron onto the parallelogram; the map
$\ygp$ maps the zero- and unit-forms of $\yOTT$ onto those of $\yOC$, while
the other four extreme forms of $\yOC$ are the images of the following
forms on $\yTT$ in the notation of fig.~\ref{fig:stat_red}:
\begin{figure}[!t]
\centering\hspace*{\fill}
\includegraphics[width=0.49\columnwidth]{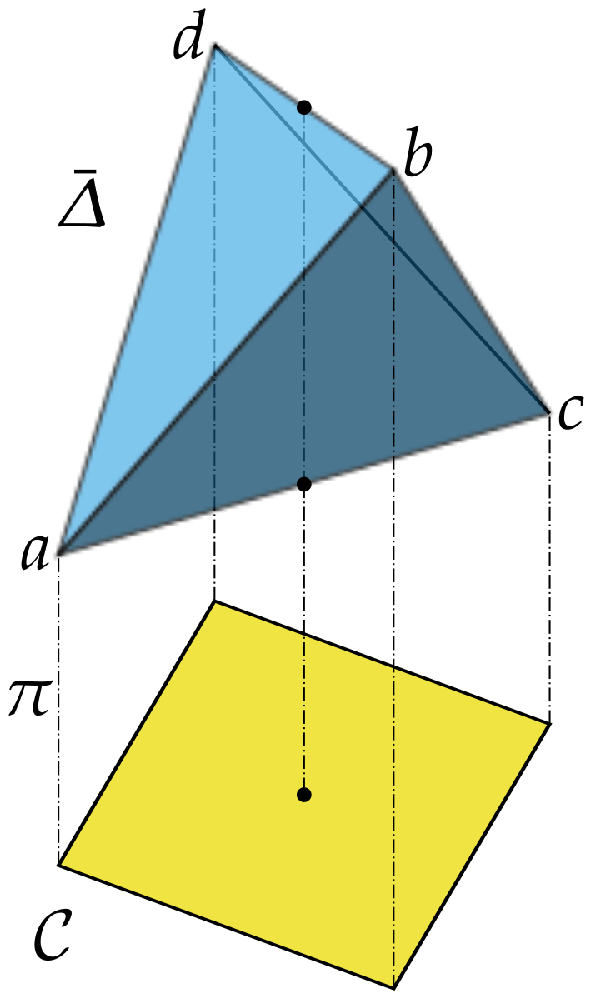}\hfill
\includegraphics[width=0.49\columnwidth]{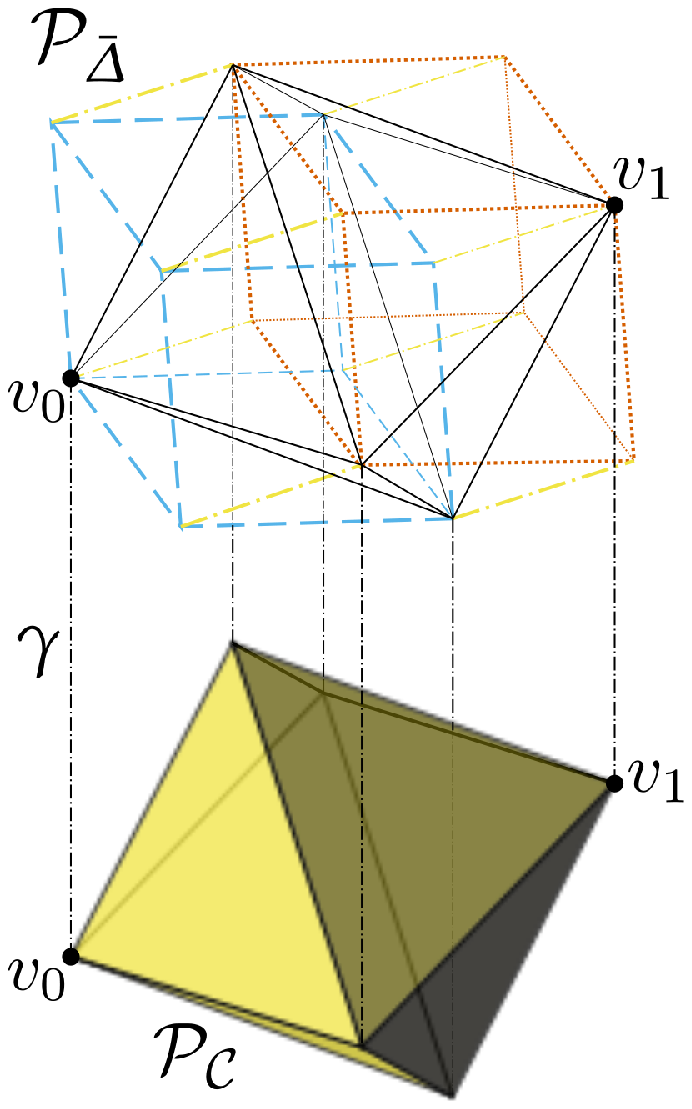}\hspace*{\fill}
\caption{\small Illustration of the refinement of the statistical model of
  Example~\ref{ex:simplic}. Note how the mapping $\yfp\colon \yTT \to \yC$
  is total and non-injective, and $\ygp\colon \yOTT \hookrightarrow \yOC$
  is partial; both are projections. The representation of the hypercube or
  four-dimensional parallelotope $\yOTT$ is explained in
  fig.~\ref{fig:4cube} on the next page.}
\label{fig:stat_red}
\end{figure}
\begin{figure}[!ph]
\centering
\includegraphics[width=\columnwidth]{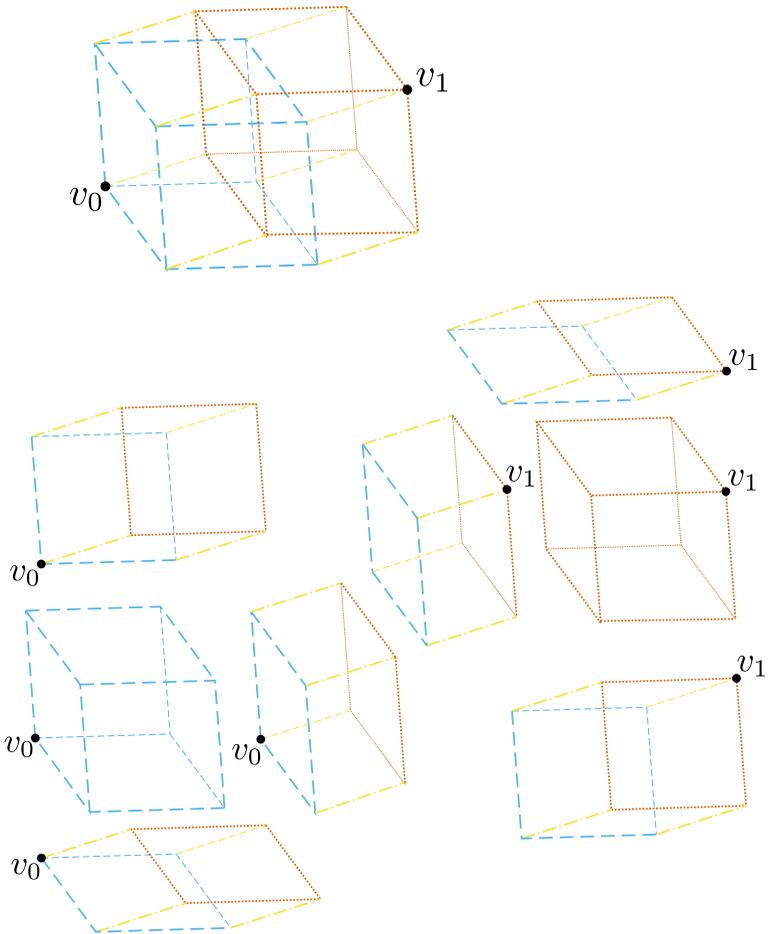}
\caption{\small A hypercube, or four-dimensional parallelotope, is a
  four-dimensional polytope with eight, pairwise parallel,
  three-dimensional facets, all parallelepipeds; and 24 two-dimensional
  faces, all parallelograms. The figure on the top is the projection of a
  4-parallelotope onto a three-dimensional space (further projected on
  paper); because of the dimensional reduction some of the projected facets
  intersect each other. To help you distinguish all eight of them in the
  top figure, they are separately represented underneath it. As the
  convex-form space of a three-dimensional simplex (tetrahedron), two
  vertices of the hypercube represent the nought- and unit-forms, also
  indicated in the figure.}
\label{fig:4cube}
\end{figure}
\begin{itemize}\tightlist
  \item that having zero value on the vertices
  $a$ and $b$ and unit value on $c$ and $d$,
\item as the previous but with zero and unit values exchanged,
\item that having zero value on the vertices
  $a$ and $d$ and unit value on $b$ and $c$,
\item as the previous but with zero and unit values exchanged.
\end{itemize}
We see that $\yfp$ is a total map, defined on all $\yTT$; whereas $\ygp$ is
partial: in particular, it is not defined on the non-constant extreme forms
of $\yOTT$.\qem
\end{example}

The preceding example is based on the fact that any convex polytope with
$m$ extreme points can be obtained as the \emph{projection} of an
$(m-1)$-dimensional simplex: see again Gr\"unbaum
\citey{gruenbaum1967_r2003}, \sect~5.1, theorem~2. We have the following
\begin{fact}
  The theorem just mentioned can always be used to construct a simplicial
  refinement, in the guise of Example~\ref{ex:simplic}, of \emph{any}
  statistical model; \cf\ Holevo \citey[\sect~I.5]{holevo1980_t1982}.
\end{fact}

\begin{figure}[!t]
\centering
\includegraphics[width=\columnwidth]{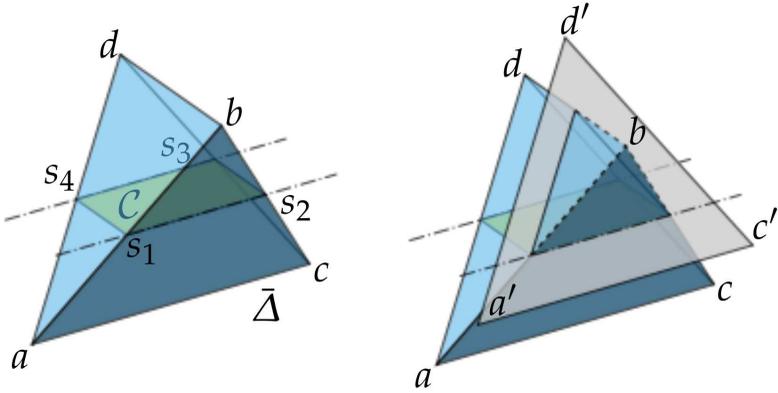}
\caption{\small An extreme convex form $v$ of the two-dimensional
  parallelogram $\yC$ is represented by the two parallel dot-dashed lines
  $s_1s_2$ and $s_3s_4$ lying in the same plane as $\yC$ (left figure). A
  convex form $w$ of $\yTT$ is represented by two parallel planes, and if
  $w$ is to correspond to $v$, $\yg(w)=v$, in such a way that $\yg(w)\inn
  s_i = v\inn s_i$, these two planes must contain $s_1s_2$ and $s_3s_4$;
  they cannot intersect the interior of $\yTT$, however, if $w$ is to be a
  convex form. But it is impossible to satisfy both requirements for both
  planes: \eg, the only plane that contains the line $s_3s_4$ and does not
  cut $\yTT$ is $acd$ (right figure); then $a'c'd'$ is the parallel plane
  containing $s_1s_2$, but this plane cuts $\yTT$ (meaning, \eg, that $w
  \inn b<0$ or $w \inn b>1$). All other constructions one can think of have
  the same problem. Thus the map $\yg$ cannot
  exist.}
\label{fig:impossquare}
\end{figure}

We already said that another theorem of convex geometry states that any
convex polytope with $m$ facets can be obtained as the \emph{section} of an
$(m-1)$-\bd dimensional simplex \cite[\sect~5.1,
theorem~1]{gruenbaum1967_r2003}. This theorem \emph{cannot} be used to
construct a simplicial refinement, though, as shown in the following
\begin{cexample}\label{ex:counter}
  The parallelogram $\yC$ of the preceding example can be obtained as the
  intersection of the tetrahedron $\yTT$ and an intersecting plane parallel
  to the segments $ac$ and $bd$ of fig.~\ref{fig:impossquare}. This defines
  a partial, surjective, affine map $\yf\colon \yTT \hookrightarrow \yC$,
  which is simply the identity, $\yf(s) =s$, in its domain of definition
  $\yC \subset \yTT$. However, it is impossible to find an affine map
  $\yg\colon \yOTT \hookrightarrow \yOC$ that be surjective and such that
  $\yom \inn s = \yg(\yom)\inn \yf(s) \equiv\yg(\yom) \inn s$: the extreme
  forms of $\yOC$, for example, cannot have any counter-image. The reason
  for this is geometrically explained in fig.~\ref{fig:impossquare}. Thus
  we cannot construct a simplicial refinement of $(\yC, \yOC)$ where $\yC$
  is obtained by sectioning $\yTT$.\qem
\end{cexample}

Linusson \citey{linusson2007} has shown that the previous counter-example
is generally valid:
\begin{theorem}[Linusson]
  Given a \emph{non-simplicial} statistical model $(\yC, \yOC)$ it is
  impossible to find a simplicial refinement $(\yT, \yOT, \yf, \yg)$ such
  that $\yf$ is injective in its domain of definition (which means that
  $\yf$ would represent the intersection of $\yT$ with a hyperplane,
  whereby $\yC$ is obtained).
\end{theorem}

\begin{example}\label{ex:pentag1}
  Let $\yC$ be a two-dimensional pentagonal convex space with extreme
  points $\set{s_1, s_2, s_3, s_4, s_5}$. Its convex-form space $\yOC$ is a
  pentagonal trapezohedron that has, besides the null- and unit-forms
  $\yoz$ and $\you$, ten other extreme points given by the forms
  \[\set{v_1,\dotsc,v_5,\you-v_1,\dotsc,\you-v_5}\] such that
\begin{equation}\label{eq:penta}
  \begin{aligned}
    (v_i \inn s_j) &=
    \begin{pmatrix}
      1&\yga&0&0&\yga \\
      \yga&1&\yga&0&0 \\
      0&\yga&1&\yga&0 \\
      0&0&\yga&1&\yga \\
      \yga&0&0&\yga&1
    \end{pmatrix} \qquad\text{with $\yga\defd \frac{\sqrt{5}-1}{2}$},
    \\
    (\you-v_i) \inn s_j &= \you \inn s_j - v_i \inn s_j= 1- v_i \inn s_j;
  \end{aligned}
\end{equation}
see fig.~\ref{fig:pentag}.
\begin{figure}[!t]
  \centering
  \includegraphics[width=\columnwidth]{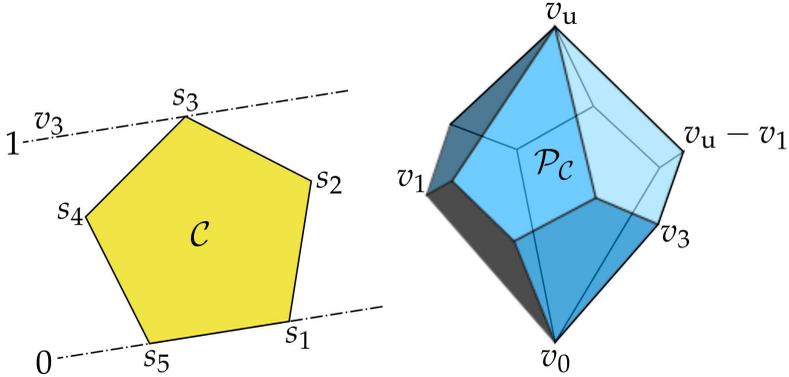}
  \caption{\small The pentagonal convex set $\yC$ and its convex-form space
    $\yOC$ from Example~\ref{ex:pentag1}. The convex form $v_3$ is shown on
    $\yC$ (as a pair of parallel lines) and on $\yOC$ (as a point).}
\label{fig:pentag}
\end{figure}

A refinement of this model is given by $(\yT, \yOT, \yf, \yg)$ where $\yT$
is a nine-dimensional simplex, or decatope, with ten extreme points
$\set{e_1, \dotsc, e_{10}}$, and $ \yOT$ is a ten-dimensional parallelotope
with twelve extreme points given by the forms \[\set{\yozz, \youu, d_1,
  \dotsc, d_{10}, \youu-d_1, \dotsc, \youu-d_{10}}\] such that
\begin{equation}
\yozz \inn e_j = 0,\quad \youu \inn e_j =1,\quad
d_i \inn e_j = \delt_{ij},\quad
(\youu-d_i)\inn e_j = 1-\delt_{ij}.
\label{eq:extr10}
\end{equation}
The map $\yf$ is defined on the points
\begin{equation}
  \label{eq:corr_pure}
  \begin{gathered}
\yf[(e_1 + e_2)/2] = s_1,\quad
\yf[(e_3 + e_4)/2] =s_2,\quad
\yf[(e_5 + e_6)/2] = s_3,\\
\yf[(e_7 + e_8)/2] = s_4,\quad
\yf[(e_9 + e_{10})/2] = s_5,
  \end{gathered}
\end{equation}
and their convex combinations; \ie, it is partial and defined on a
four-dimensional simplex given by the convex span of $\set{(e_1 + e_2)/2,
  (e_3 + e_4)/2,\dotsc,(e_9 + e_{10})/2}$; see fig.~\ref{fig:pentagproj}.
\begin{figure}[!t]
  \centering
  \includegraphics[width=\columnwidth]{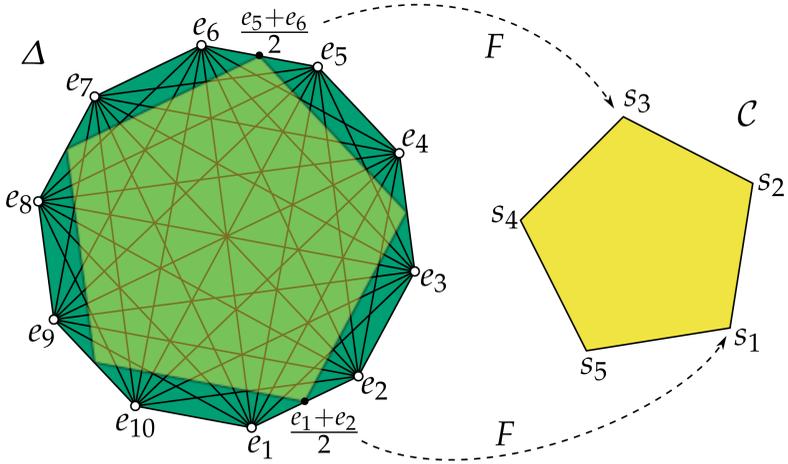}
  \caption{\small Representation of the decatope $\yT$ and the map $\yf$ of
    Example~\ref{ex:pentag1}. $\yT$ is represented by its graph
    \citep[\sects~11.3, 8.4]{gruenbaum1967_r2003}, which can be understood
    as a parallel projection of $\yT$ onto a two-dimensional plane: the
    vertices represent its extreme points $\set{e_i}$, the
    $\tbinom{10}{2}=45$ lines connecting any two vertices represent its
    edges, and all $\tbinom{10}{3}=120$ triangles connecting any three
    vertices represent its two-dimensional faces. The map $\yf$ is only
    defined on a four-dimensional, simplicial convex subset (in lighter
    green) of $\yT$.}
\label{fig:pentagproj}
\end{figure}

The map $\yg$ is also partial:
define the following one-dimensional convex subsets
\begin{equation}
  \label{eq:corr_outc}
  \begin{gathered}
    \begin{aligned}
      \ydr_1 &\defd \set{d_1 + d_2 + \yb (d_9 + d_3)+ \yc (d_{10} + d_4)\st \yb \in \clcl{2\yga -1,1}},\\
      \ydr_2 &\defd \set{d_3 + d_4 + \yb (d_1 + d_5)+ \yc (d_2 + d_6)\st \yb \in \clcl{2\yga -1,1}},\\
      \ydr_3 &\defd \set{d_5 + d_6 + \yb (d_3 + d_7)+ \yc (d_4 + d_8)\st \yb \in \clcl{2\yga -1,1}}, \\
      \ydr_4 &\defd \set{d_7 + d_8 + \yb (d_5 + d_9)+ \yc (d_6 + d_{10})\st \yb \in \clcl{2\yga -1,1}}, \\
      \ydr_5 &\defd \set{d_9 + d_{10} + \yb (d_7 + d_1)+ \yc (d_8 +
        d_2)\st \yb \in \clcl{2\yga -1,1}},
    \end{aligned}
    \\
    \text{with $\yc\defd 2 \yga -\yb$};
  \end{gathered}
\end{equation}
then $\yg$ is defined by
\begin{equation}
  \label{eq:corr_outc}
  \begin{gathered}
    \yg(\yozz) = \yoz,\qquad \yg(\youu) = \you,\qquad \yg(\ydr_i) =
    \set{v_i},\; i=1,\dotsc,5, \\ \yg(\youu - \ydr_i) =
    \set{\you-v_i},\; i=1,\dotsc,5,
  \end{gathered}
\end{equation}
where $\youu-\ydr_i \defd \set{\youu-a \st a\in \ydr_i}$.

The important features of this refinement are these:
\begin{enumerate}[a.]
\item each extreme point of $\yC$ corresponds to a \emph{non-extreme} point
  of $\yT$, and to that alone;
\item each of the ten non-trivial extreme points of $\yOC$ corresponds to
  a non-zero-dimensional convex set of non-extreme points of $\yOT$;
\item the map $\yf$ is partial, \ie\ it is not defined on some points of
  $\yT$, not even its ten pure ones $\set{e_i}$; contrast this with
  Example~\ref{ex:simplic};
\item on the four-dimensional simplex on which it is defined, the map $\yf$
  acts as a projection onto $\yC$, analogously to the map $\yfp$ of
  Example~\ref{ex:simplic}.\qem
\end{enumerate}
\end{example}

In the example just discussed, the fact that the extreme points of $\yC$
correspond to single points of $\yT$ makes it possible for the extreme
points of $\yOC$ to correspond to one-dimensional sets in $\yOT$. But the
opposite situation is also possible:

\begin{example}\label{ex:pentag2}
The sets $\yC$, $\yOC$, $\yT$, $\yOT$ are defined as in the preceding
example, but the maps $\yf$ and $\yg$ are defined differently. Consider these
five one-dimensional faces of $\yT$:
\begin{equation}
  \label{eq:corr_outc}
  \begin{gathered}
\yer_1 \defd \conv\set{e_1,e_2} = \set{\yb e_1 + (1-\yb)e_2 \st x\in\clcl{0,1}}\\
\begin{aligned}
  \yer_2 &\defd \conv\set{e_3,e_4}, \qquad&
  \yer_3 &\defd \conv\set{e_5,e_6},\\
  \yer_4 &\defd \conv\set{e_7,e_8},\qquad &\yer_5 &\defd
  \conv\set{e_9,e_{10}};
\end{aligned}
\end{gathered}
\end{equation}
then define $\yf$ by
\begin{equation}
  \label{eq:corr_outc}
\yf(\yer_i) = \set{s_i},\; i=1,\dotsc,5,
\end{equation}
see fig.~\ref{fig:pentag2},
\begin{figure}[!b]
  \centering
  \includegraphics[width=\columnwidth]{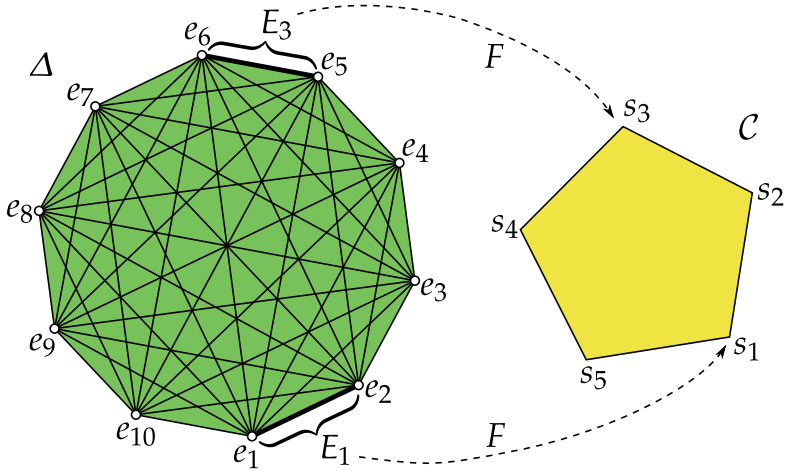}
  \caption{\small Representation of the decatope $\yT$ and the map $\yf$ of
    Example~\ref{ex:pentag2}; \cf\ fig.~\ref{fig:pentagproj}. The map $\yf$
    now maps five edges of $\yT$ to the $\set{s_i}$ and by convex
    combination is defined on all of $\yT$.}
\label{fig:pentag2}
\end{figure}
and $\yg$ by
\begin{equation}
  \label{eq:corr_outc_fix}
  \begin{gathered}
    \yg(\yozz) = \yoz,\quad
    \yg(\youu) = \you,\\
    \begin{aligned}
      &\yg[d_1 + d_2 + \yga (d_3 + d_9)+ \yga (d_4 + d_{10})] = r_1,\\
      &\yg[   d_3 + d_4 + \yga (d_1 + d_5)+ \yga (d_2 + d_6)]= r_2,\\
      &\yg[  d_5 + d_6 + \yga (d_3 + d_7)+ \yga (d_4 + d_8)]=r_3,\\
      &\yg[  d_7 + d_8 + \yga (d_5 + d_9)+ \yga (d_6 + d_{10})]=r_4,\\
      &\yg[ d_9 + d_{10} + \yga (d_1 + d_7)+ \yga (d_2 + d_8)]=r_5.
    \end{aligned}
  \end{gathered}
\end{equation}

This refinement differs from the one in the previous Example in that
\begin{enumerate}[a.]
\item the five extreme points of $\yC$ correspond to five
  one-dimensional faces of $\yT$;
\item the ten non-trivial extreme points of $\yOC$ correspond to single,
  non-extreme points of $\yOT$;
\item the map $\yf$ is total, it is indeed a parallel projection of $\yT$
  onto $\yC$.
\end{enumerate}
This refinement is in fact more similar to that of
Example~\ref{ex:simplic}, with the difference that the simplex from which
the polytope is obtained by projection is not the one of the least possible
dimension. \qem
\end{example}

\chapter{Conjectures and questions}
\label{sec:conjectures}

I do not know of any study of the general properties of simplicial
refinements of a statistical model. Apart from Linusson's theorem, general
theorems are lacking.

It seems that if we try to construct a refinement of a statistical model
$(\yC, \yOC)$ with a simplex $\yT$ having fewer extreme points than $\yC$,
the construction runs into problems similar to those of Counter-\bd
example~\ref{ex:counter} for the map $\yg$. Moreover, even considering
higher-\bd dimensional simplices, it seems that if we try to construct
$\yf\colon \yT\hookrightarrow \yC$ in such a way that only some
\emph{non-extreme} points of $\yT$ map onto some extreme points of $\yC$,
as in Example~\ref{ex:pentag1}, then those non-extreme points have to be
enough `far apart' face-wise, otherwise we cannot construct $\yg$, again
for problems like those in Counter-example~\ref{ex:counter}.

These remarks naturally lead to the following conjectures. Unfortunately I
have not been able to prove or disprove them, and wish that convex
geometers will take note of them:

\bigskip
Consider a non-simplicial statistical model $(\yC, \yOC)$,
where $\yC$ has $m$ extreme points:
 \begin{conj}\label{con:convconj}
   All simplicial refinements of $(\yC, \yOC)$ have the form $(\yT, \yOT,
   \yff\circ \yfp, \ygg\circ\ygp)$, where:
\begin{enumerate}[a.]
\item $(\yTT, \yOTT, \yfp, \ygp)$ is the refinement of $(\yC, \yOC)$ where
  $\yTT$ is the $(m-1)$-dimensional simplex from which $\yC$ is obtained by
  parallel projection $\yfp$, as in Example~\ref{ex:simplic},
\item $(\yT, \yOT, \yff, \ygg)$ is a refinement of $(\yT, \yOT)$, where
  $\yT$ is a simplex of dimension larger than $(m-1)$ and $\yff\colon \yT
  \hookrightarrow \yTT$ is an affine, onto, partial map (a partial
  projection, an intersection, or a combination of the two).
\end{enumerate}
In other words, the refinement obtained by projection of an
$(m-1)$-dimensional simplex is the lowest-dimensional one. See
fig.~\ref{fig:convconj}.
\begin{figure}[!b]
  \centering
  \includegraphics[width=\columnwidth]{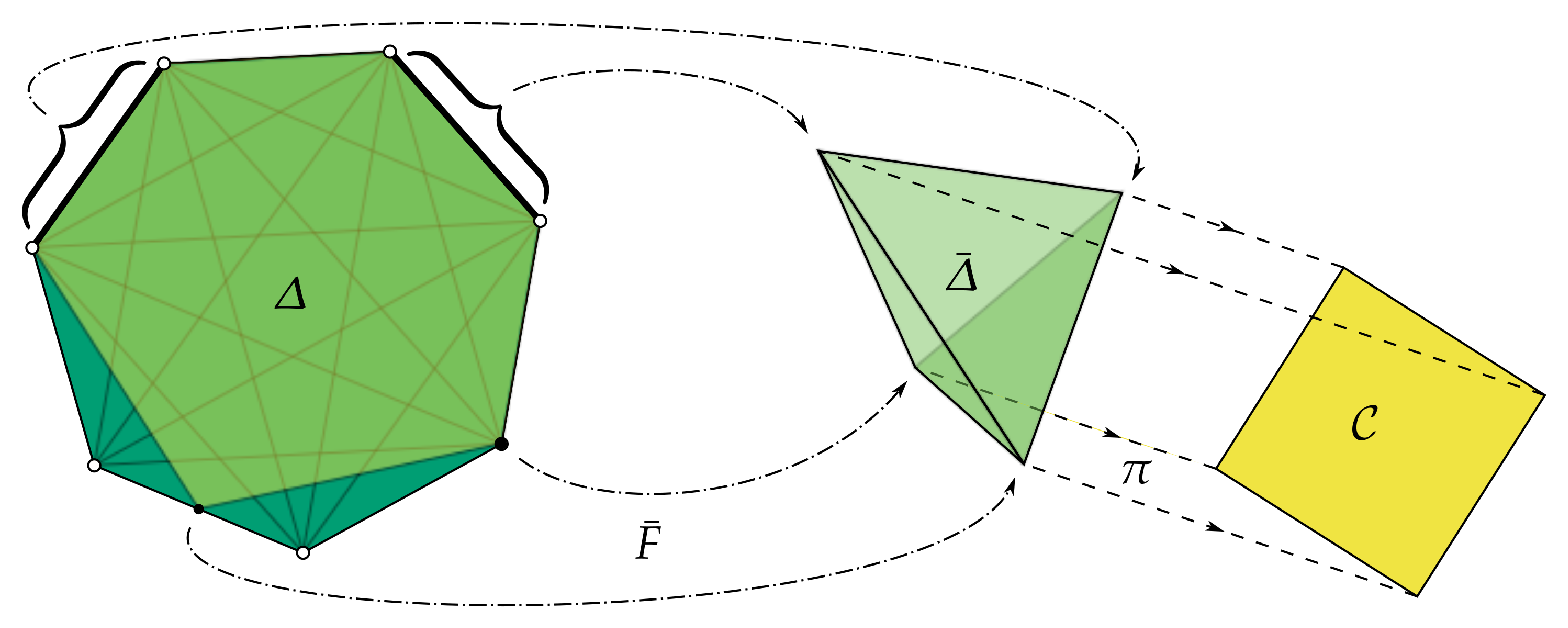}
  \caption{\small Conjecture~\ref{con:convconj} says that any refinement
    $(\yT, \yOT, \yf, \yg)$ of $(\yC, \yOC)$ has the schema of the above
    figure, with $\yf=\yff\circ\yfp$ and $\yg=\ygg\circ \ygp$. In other words,
    the refinement obtained by projection is the one of lowest dimension.}
\label{fig:convconj}
\end{figure}
 \end{conj}
 \begin{conj}\label{conj:ep}
   No simplicial refinement $(\yT, \yOT, \yf, \yg)$ exists with $\yT$
   having fewer extreme points than $\yC$. This is also a corollary of the
   previous conjecture.
 \end{conj}
 \begin{conj}\label{conj:faces}
   Let $(\yT, \yOT, \yf, \yg)$ be any refinement of $(\yC, \yOC)$, and
   $a,b$ be any two extreme points of $\yC$. Let $\yf^{-1}({a}),
   \yf^{-1}({b})$ be their counter-images in $\yT$, and $A,B$ the minimal
   faces of $\yT$ containing these counter-images. Then $A\cap B
   =\emptyset$. In other words, no refinement exists such that two extreme
   points of $\yC$ correspond to points in $\yT$ lying on adjacent faces.
 \end{conj}


 \bigskip The following are very important questions for the application of
 the theory of statistical models to quantum theory and general physical
 statistical theories. Given a non-simplicial statistical model $(\yC,
 \yOC)$:
\begin{quest}
  Does a simplicial refinement exist with $\yT$ having fewer extreme points
  than $\yC$? (\Cf\ Conjecture~\ref{conj:ep}.)
\end{quest}
\begin{quest}
  What is the least dimension of $\yT$ among all simplicial refinements?
\end{quest}
\begin{quest}
  For each simplicial refinement, can the affine partial map $\yf$ be
  extended to a map $\aff\yT \to \aff\yC$? If so, is the extension unique?
  What about an analogous extension of $\yg$?
\end{quest}
\begin{quest}
  How to extend the theory of simplicial refinements to statistical models
  $(\yC,\yOC)$ where $\yC$ has a continuum of extreme points or even
  infinite dimension?
\end{quest}

\appendix

\chapter*{Appendix: Physical motivations}
\label{sec:phys}

In quantum theory and other physical statistical theories the set of
statistical states of a system has the structure of a convex set $\yC$
(generally having a continuum of extreme points, and in some cases infinite
dimension). By statistical states I mean, \eg, the statistical operators in
quantum theory or the Liouville distributions in statistical mechanics. The
set of measurement outcomes, for all kinds of measurement that can be made
on the system, has the convex structure of the set $\yOC$ (strictly
speaking this only holds for \emph{maximal} theories
\citep{portamana2011c}, but all known physical theories are maximal). The
probability of obtaining the outcome represented by a form $v \in \yOC$
given that the measurement associated to the same form is made on the state
represented by a point $s \in \yC$ is then given by the action of the form
$v$ on the point $s$, $v\inn s$. The probabilistic features of a physical
system are thus represented by a statistical model $(\yC, \yOC)$. In what
follows I assume these notions valid (with topological care) for convex
spaces with a continuum of extreme points and possibly infinite dimensions.

The convex structure of the state space $\yC$ of quantum and other
statistical theories is the origin of many peculiar statistical features,
\eg\ the fact that some pure states (\ie\ states represented by extreme
points) cannot be distinguished by one measurement instance, or the fact
that no measurement can resolve situations of uncertainties between two or
more particular, different sets of pure states. Most classical systems do
not present these peculiarities: because their state spaces have simplicial
convex structures $\yT$.

The connexions between states and measurements of two physical systems are
translated into the connexions between the statistical models associated to
them. This is a case that has interested physicists for many years: to
reproduce the probabilistic features of quantum systems as emergent from
those of classical ones; more or less as it happens with thermostatics and
statistical mechanics. The quantum system's phenomenology would then be a
`coarser' version of the classical system's one.

This relationship between a `coarser' physical system represented by $(\yC,
\yOC)$ and a `finer' classical one represented by $(\yT, \yOT)$ implies a
correspondence between their states and measurement outcomes with several
natural requirements:
\begin{enumerate}[a.]
\item all states in the coarser system must have a correspondent in the
  finer, otherwise the latter system would not be able to describe all the
  phenomenology of the former;
\item the finer system can have states with no counterpart in the coarser,
  because its phenomenology can be
  richer;
\item several states in the finer system can correspond to the same state
  in the coarser: in the coarser system we simply cannot tell them apart;
\item at most one state in the finer system may correspond to one in the
  coarser, because the finer system cannot have less `resolving power' than
  the coarser;
\item a situation of uncertainty among states in the finer system must
  correspond, probability-wise, to the same situation of uncertainty among
  the corresponding states in the coarser, otherwise the finer system could
  not reproduce the statistical phenomenology of the coarser;
\end{enumerate}
and analogous requirements hold for the measurement outcomes of the two
systems. Moreover, the outcome probabilities for all states and measurement
outcomes of the coarser system must obviously be reproduced by the finer
one if the statistical phenomenology of the former is to be obtained. We
immediately see that these requirements are translated into the
characteristics C1--C5 and \eqn~\eqref{eq:compatib} of \sect~\ref{sec:maps}
for the statistical models associated to the coarser and finer systems. In
other words, we are saying that the statistical model associated to the
finer classical system is a \emph{simplicial refinement} $(\yT, \yOT,
\yf,\yg)$ of that associated to the coarser physical system $(\yC, \yOC)$.

The best-known example of a classical-like theory
reproducing quantum systems described by Schr\"odinger's equation is
Bohmian mechanics \citep{bohm1952a,bohm1952b,berndletal1995}. Its pure
states are of the form $(\psi, \yx)$, where the $\set{\psi}$ are in one-one
correspondence with quantum pure states (wave functions) and the
$\set{\yx}$ are extra configuration variables. A quantum pure state $\psi$
corresponds in Bohm's theory to the mixed state given by the distribution
\begin{equation}
\yll(\psi',\yx')\,\di\psi'\di\yx'
\defd \delt(\psi'-\psi) \abs{\psi(\yx')}^2\,\di\psi'\di\yx';
\label{eq:bohmq}
\end{equation}
thus the refinement given by Bohm's theory behaves analogously to
Example~\ref{ex:pentag1} of \sect~\ref{sec:maps}.

The conjectures of the last section, extended to the case of convex sets
with a continuum of extreme points and possibly infinite dimensions,
translate into physical requirements for a classical system to be a finer
version of generic statistical one, a quantum one in particular:
\begin{tconj}
  Given a generic physical statistical system, \eg\ a quantum one, the
  simplest classical physical system capable of reproducing its statistical
  phenomenology is one whose set of pure states can be put in one-one
  correspondence with the set of pure states of the generic one. In the
  quantum case this says that models like the Beltrametti-Bugajski one
  \citey{beltramettietal1995}, which use Holevo's construction
  \citey[\sect~I.5]{holevo1980_t1982}, are the simplest possible.
\end{tconj}

\begin{tconj}
  If a classical system reproduces the statistical phenomenology of a
  generic physical statistical, \eg\ a quantum one, then its set of pure
  states cannot have lower cardinality that the generic system's one. In
  the quantum case, the manifold of pure states of the classical system
  cannot have less dimensions than that of the quantum one. This means that
  behind the phenomenology of Schr\"odinger's equation there must be
  \emph{physical fields}: particles only are not enough.
\end{tconj}
This is true of all classical models that reproduce quantum ones. For
example, the set of pure states of Bohmian mechanics obviously has a
higher cardinality than the quantum one. Also other models, like Nelson's
\citey{nelson1966,nelson1985}, have a greater cardinality than the quantum
systems they reproduce, because beside the configuration variables $\yx$
they introduce other quantities that can only be interpreted as physical
fields from the way they enter into the equations of motion.

\begin{tconj}
  Let $\ypa, \ypb$ be two pure states of the generic statistical system
  (\eg, again, a quantum one), and $\yf^{-1}({\ypa}), \yf^{-1}({\ypb})$ the
  sets of statistical classical states reproducing those states. Let $A,B$
  be the sets of pure classical states into which $\yf^{-1}({\ypa}),
  \yf^{-1}({\ypb})$ can be purified. Then the sets $A,B$ have neither pure
  nor mixed states in common.
\end{tconj}
This conjecture is also satisfied in Nelson's model or in Bohmian
mechanics, owing to the delta function in \eqn~\eqref{eq:bohmq} that makes
two distributions corresponding to two different quantum pure states to
have disjoint supports. Proving Conjecture~\ref{conj:faces} would rule out
models that Harrigan \amp\ Spekkens \citey{harriganetal2007_r2010} call
`$\psi$-epistemic', not only for quantum theory but for \emph{any}
non-classical (\ie, non-simplicial) maximal physical statistical theory.

\begin{acknowledgements}
  Many thanks to the developers and maintainers of \LaTeX, Emacs, AUC\TeX,
  MiK\TeX, arXiv, Inkscape; to Gen Kimura, Ingemar Bengtsson, Lucien Hardy,
  Mariano Cadoni for their interest in this study; and especially to Svante
  Linusson for his contributions to it during its first stages. My hate
  goes to
  \rotatebox{15}{\reflectbox{P}}\rotatebox{5}{I}\rotatebox{-10}{P}\rotatebox{10}{\reflectbox{P}}\rotatebox{-5}{O}.
  Research at the Perimeter Institute is supported by the Government of
  Canada through Industry Canada and by the Province of Ontario through the
  Ministry of Research and Innovation.
\end{acknowledgements}


\pagebreak[1]
\defbibnote{prenote}{%
\small
}

\newcommand{\citein}[2][]{\textnormal{\textcite[#1]{#2}}\addtocategory{1}{#2}}
\newcommand{\citebi}[2][]{\textcite[#1]{#2}\addtocategory{1}{#2}}
\defbibfilter{1}{\category{1} \or \segment{1}}
\printbibliography[filter=1,prenote=prenote]


\end{document}